\documentclass[letterpaper, 10 pt, conference]{ieeeconf}  

\IEEEoverridecommandlockouts                              
\usepackage{graphics} 
\usepackage{epsfig} 
\usepackage{mathptmx} 
\usepackage{times} 
\usepackage{amsmath} 
\usepackage{amssymb}  
\usepackage{color}
\usepackage{graphicx,subfigure}

\newtheorem{remark}{Remark}[section]

\title{\LARGE \bf A new model-free design for vehicle control\\ and its
validation through an advanced simulation platform}
 
\author{Lghani MENHOUR$^\text{a}$, Brigitte d'ANDR\'EA-NOVEL$^\text{b}$, Michel FLIESS$^\text{c, d}$, Dominique GRUYER$^\text{e}$ ~\\ 
and Hugues MOUNIER$^\text{f}$
\thanks{This work was partially supported by the French national project INOVE/ANR 2010 BLANC 308.
}
\thanks{$^\text{a}$CReSTIC, IUT de Troyes, 9 rue du Qu\'ebec, 10026 Troyes, France. \newline {\tt\small lghani.menhour@univ-reims.fr}
}
\thanks{$^\text{b}$Centre de Robotique, Mines ParisTech, PSL Research University, 60 boulevard Saint-Michel, 75272 Paris cedex 06, France. \newline {\tt\small brigitte.dandrea-novel@mines-paristech.fr}
}
\thanks{$^\text{c}$LIX (CNRS, UMR 7161), \'Ecole polytechnique, 91128 Palaiseau, France. \newline {\tt\small Michel.Fliess@polytechnique.edu}
}
\thanks{$^\text{d}$AL.I.E.N. (ALg\`{e}bre pour Identification \& Estimation Num\'{e}riques), 24-30 rue Lionnois, BP 60120, 54003 Nancy, France. \newline {\tt\small michel.fliess@alien-sas.com}
}
\thanks{$^\text{e}$IFSTTAR-CoSys-LIVIC, 77 rue des Chantiers, 78000 Versailles, France. \newline {\tt\small dominique.gruyer@ifsttar.fr}
}
\thanks{$^\text{f}$L2S (UMR 8506), CNRS -- Sup\'elec -- Universit\'{e} Paris-Sud, 3 rue Joliot-Curie, 91192 Gif-sur-Yvette, France. \newline {\tt\small hugues.mounier@lss.supelec.fr}
} 
}

\begin{document}

\maketitle
\thispagestyle{empty}
\pagestyle{empty}

\begin{abstract}
A new model-free setting and the corresponding ``intelligent'' P and PD controllers are employed for the longitudinal and lateral motions of a vehicle. This new approach has been developed and used in order to ensure simultaneously a best profile tracking for the longitudinal and lateral behaviors. The longitudinal speed and the derivative of the lateral deviation, on one hand, the driving/braking torque and the steering angle, on the other hand, are respectively the output and the input variables. Let us emphasize that a ``good'' mathematical modeling, which is quite difficult, if not impossible to obtain, is not needed for such a design. An important part of this publication is focused on the presentation of simulation results with actual and virtual data. The actual data, used in Matlab as reference trajectories, have been obtained from a properly instrumented car (Peugeot 406). Other virtual sets of data have been generated through the interconnected platform SiVIC/RTMaps. It is a dedicated virtual simulation platform for prototyping and validation of advanced driving assistance systems.


\textbf{\textsl{Keywords}}---{\hspace {0.05cm}} Longitudinal and lateral vehicle control, model-free control, intelligent P controller (i-P controller), algebraic estimation, ADAS (Advanced Driving Assistance Systems).

\end{abstract}

\section{Introduction}\label{introd}

The control of the longitudinal and lateral motions of a vehicle has been widely investigated via model-based techniques (see, \textit{e.g.}, \cite{Khodayari10,Manceur13, Menhour13a, Nobe01, Odenthal99, Poussot11a, Villagra09}, and the references therein). It should be however emphasized that obtaining a ``good'' mathematical modeling is a difficult task, if not an impossible one, since complex uncertainties and disturbances, like, for instance, frictions, should be taken into account. This is why a first model-free setting \cite{ijc13} was recently proposed in \cite{Menhour13b}. Although the corresponding numerical results were good, this attempt was suffering from the fact that one of the two flat outputs was depending on some uncertain parameters and from the fact that an accurate tracking trajectory is not guaranteed, especially when the trajectories are characterized by tight bends. 

Those tedious features vanish here thanks to a new model-free approach where: 
\begin{itemize}
\item the need to exploit the flatness property of a simplified model disappears,
\item the flat output depending on uncertain parameters is replaced by a simpler and perhaps, more natural lateral deviation.
\end{itemize}

Our paper is organized as follows. Section \ref{Section_1} gives a short summary on the model-free control approach. It is exploited in Section \ref{Section_4} in order to describe a new longitudinal and lateral vehicle control. 
Section \ref{Section_5} displays several numerical simulations, first by using real data from a 10DoF properly instrumented Peugeot 406 vehicle model, then thanks to the interconnected platforms SiVIC and RTMaps. 
This real-time simulator provides an efficient tool to replace an actual system by a virtual one. It allows advanced prototyping and validation of the perception and control algorithms with quite realistic models of the vehicle, the sensors and the environment. Some concluding remarks may be found in Section \ref{Section_6}.

\vspace{4mm}
\noindent In this paper, the following variables and notations, with their meaning will be used:
\begin{itemize}
\item $V_x$, $V_y$: longitudinal and lateral speeds in $[m/s]$, 
\item $\psi$: yaw angle in $[rad]$, 
\item $\dot{\psi}$: yaw rate in $[rad/s]$, 
\item $T_{\omega}$: acceleration/braking torque in [$Nm$], 
\item $ \delta$: steering wheel angle in [$deg$], 
\item $L_f$: distances from the CoG to the front axle in $[m]$, 
\item $L_r$: distances from the CoG to the rear axle in $[m]$, 
\item $I_z$: yaw moment of inertia in [$Kg.m^{-2}$], 
\item $g$: acceleration due to gravity $[m/s^2]$, 
\item $m$: vehicle mass in $[kg m^2]$.
\end{itemize}

\section{A short summary of model-free control\protect\footnote{See \cite{ijc13} for more details.}}
\label{Section_1}
Model-free control was already applied and used quite successfully in a lot of various concrete examples (see the references in \cite{ijc13,ecc}). For obvious reasons let us insist here on its applications to \emph{intelligent 
transportation systems}: see \cite{Abouaissa12,Choi09,Villagra09,vil2},  and \cite{Menhour13b}. This last reference was briefly discussed in Section \ref{introd}.

\subsection{The ultra-local model}\label{A}

Replace the unknown SISO system by the \emph{ultra-local model}:
\begin{equation}
y^{(\nu)} = F + \alpha u
\label{ultralocal}
\end{equation}
where 
\begin{itemize}
\item $\nu \geq 1$ is the derivation order,
\item $\alpha \in \mathbb{R}$ is chosen such that $\alpha u$ and $y^{(\nu)}$ are of the same order of magnitude, 
\item $\nu$ and $\alpha$ are chosen by the practitioner. 
\end{itemize}
\begin{remark}
In all the existing concrete examples 
$$\nu = \ 1 \ \text{or} \ 2$$ 
Until now from our knowledge, in the context of model-free control, the example of magnetic bearings \cite{Miras13} with their low friction, provides the only instance where the order $\nu = 2$ is necessary. 
\end{remark}
Some comments and assumptions on $F$ can be done:
\begin{itemize}
\item $F$ is estimated via the measurements of the control input $u$ and the controlled output $y$,
\item $F$ does not distinguish between the unknown model of the system and the perturbations and uncertainties.
\end{itemize}

\subsection{Intelligent controllers}
Set $\nu = 2$ in Equation \eqref{ultralocal}:
\begin{equation}
\label{MFC_4_n_2}
\ddot{y} = F+ \alpha u 
\end{equation}
The corresponding \emph{intelligent Proportional-Integral-Derivative controller}, or \emph{iPID}, reads
\begin{equation}
\label{iPID_c}
u = - \frac{\left( F - \ddot{y}^d + K_P e + K_I\int e dt + K_D \dot{e}\right)}{\alpha}  
\end{equation}
where
\begin{itemize}	
\item $y^d$ is the reference trajectory,
\item $e = y-y^d$ is the tracking error and $y^d$ is a desired signal, 
\item $K_P$, $K_I$, $K_D \in \mathbb{R}$ are the usual gains.
\end{itemize}
Combining Equations \eqref{MFC_4_n_2} and \eqref{iPID_c} yields
$$
\dot{e} + K_P e + K_I \int e dt + K_D \dot{e}= 0
$$
where $F$ does not appear anymore. Gain tuning becomes therefore quite straightforward. This is a major benefit when compared to ``classic'' PIDs. 
If $K_I = 0$ we obtain the \emph{intelligent Proportional-Derivative controller}, or \emph{iPD},
\begin{equation}
\label{iPD_c}
u = - \frac{\left( F - \ddot{y}^d + K_P e dt + K_D \dot{e}\right)}{\alpha}  
\end{equation}
Set $\nu = 1$ in Equation \eqref{ultralocal}:
\begin{equation}
\label{MFC_4_n_1}
\dot{y} = F+ \alpha u 
\end{equation}
The corresponding \emph{intelligent Proportional-Integral controller}, or \emph{iPI}, reads: 
\begin{equation}
\label{iPI_c}
u = - \frac{\left(F - \dot{y}^d + K_P e + K_I\int e dt \right)}{\alpha}
\end{equation}	
If $K_I = 0$ in Equation \eqref{iPI_c}, we obtain the \emph{intelligent proportional controller}, or \emph{iP}, which turns out until now to be the most useful intelligent controller:
\begin{equation}
\label{ip}
{u = - \frac{F - \dot{y}^\ast + K_P e}{\alpha}}
\end{equation} 

\subsection{Algebraic estimation of $F$}
\label{F}
$F$ in Equation \eqref{ultralocal} is assumed to be ``well'' approximated by a piecewise constant function $F_{\text{est}} $. According to the algebraic parameter identification developed in \cite{sira1,sira2}, where the probabilistic properties of the corrupting noises may be ignored, if $\nu = 1$, Equation \eqref{MFC_4_n_1} rewrites in the operational domain (see, \textit{e.g.}, \cite{Yosida84}) 
$$
s Y = \frac{\Phi}{s}+\alpha U +y(0)
$$
where $\Phi$ is a constant. We get rid of the initial condition $y(0)$ by multiplying the both sides on the left by $\frac{d}{ds}$:
$$
Y + s\frac{dY}{ds}=-\frac{\Phi}{s^2}+\alpha \frac{dU}{ds}
$$
Noise attenuation is achieved by multiplying both sides on the left by $s^{-2}$. It yields in the time domain the realtime estimation
{\small \begin{equation}\label{integral}
\small F_{\text{est}}(t)  =-\frac{6}{\tau^3}\int_{t-\tau}^t \left\lbrack (\tau -2\sigma)y(\sigma)+\alpha\sigma(\tau -\sigma)u(\sigma) \right\rbrack d\sigma 
\end{equation}}
where $\tau > 0$ might be quite small. This integral may, of course, be replaced in practice by a classic digital filter. See \cite{nice} for a cheap and small hardware implementation of our controller. The extension to the case $\nu = 2$ is straightforward.


\section{Application to a vehicle control}
\label{Section_4}
Select, in order to avoid any modeling problem, the following input and output variables: 
\begin{enumerate}
\item the acceleration/braking torque $u_1 = T_\omega$ and the longitudinal speed $y_1$,
\item the steering wheel angle  $u_2 = \delta$ and the lateral deviation $y_2$.
\end{enumerate}

\noindent Newton's second law yields then the two local models:
 \begin{eqnarray}
\label{MFC_LLVC2_1}
\text{longitudinal local model:} & \dot{y}_1 = F_1+ \alpha_1 u_1\\
\label{MFC_LLVC2_2}
\text{lateral local model:} & \ddot{y}_2 = F_2+ \alpha_2 u_2 
\end{eqnarray}
Note the following properties:
\begin{itemize}
\item Equations \eqref{MFC_LLVC2_1}-\eqref{MFC_LLVC2_2} seem decoupled, but the coupling effects are included in the terms $F_1$ and $F_2$.
\item Equation \eqref{MFC_LLVC2_2} is an order $2$ formula with respect to the derivative of $y_2$.
\end{itemize}
For Equation \eqref{MFC_LLVC2_1} (resp. \eqref{MFC_LLVC2_2}), the loop is closed by an iP \eqref{ip} (resp. iPD \eqref{iPD_c}).

\section{Simulations with actual data}
\label{Section_5}
The simulations are carried out firstly with Matlab using a 10DoF instrumented Peugeot 406 car as in \cite{Menhour13b}, and secondly thanks to the SiVIC\footnote{SiVIC is a professional software of CIVITEC \newline ({\tt http://www.civitec.com}).} simulator \cite{Gruyer10b, Gruyer10a, Vanholme10, Gruyer09} interfaced with RTMaps platform\footnote{RT-Maps is developed by Intempora \newline ({\tt http://www.intempora.com}).}.

\subsection{Simulation under Matlab}
\begin{figure}[!ht]
\centering
\includegraphics[scale=0.56]{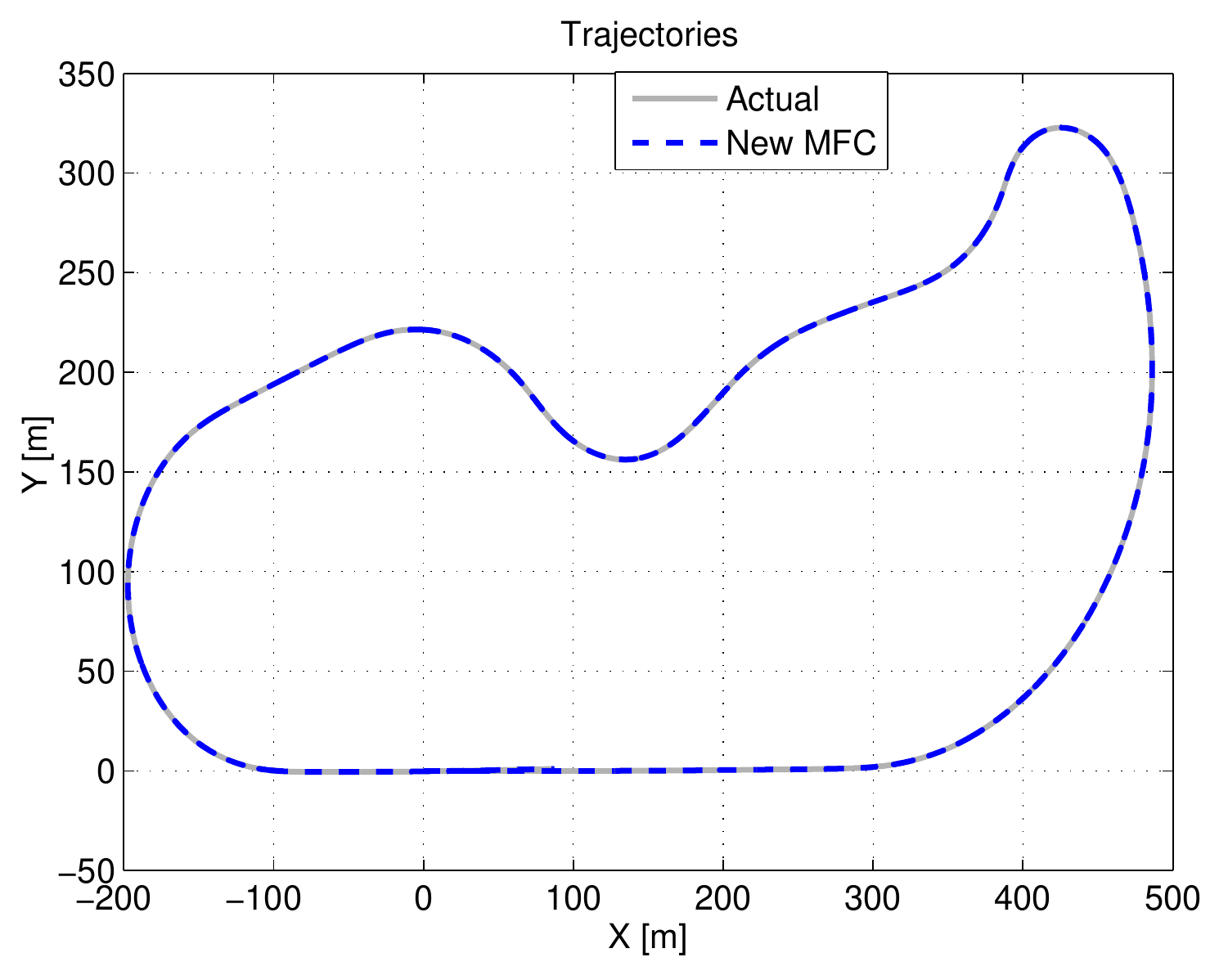}
\vspace{-0.6cm}
\caption{Reference trajectory versus the closed-loop simulated trajectory}
\label{X_Y}
\end{figure}

\begin{figure}[!ht]
\centering
\includegraphics[scale=0.56]{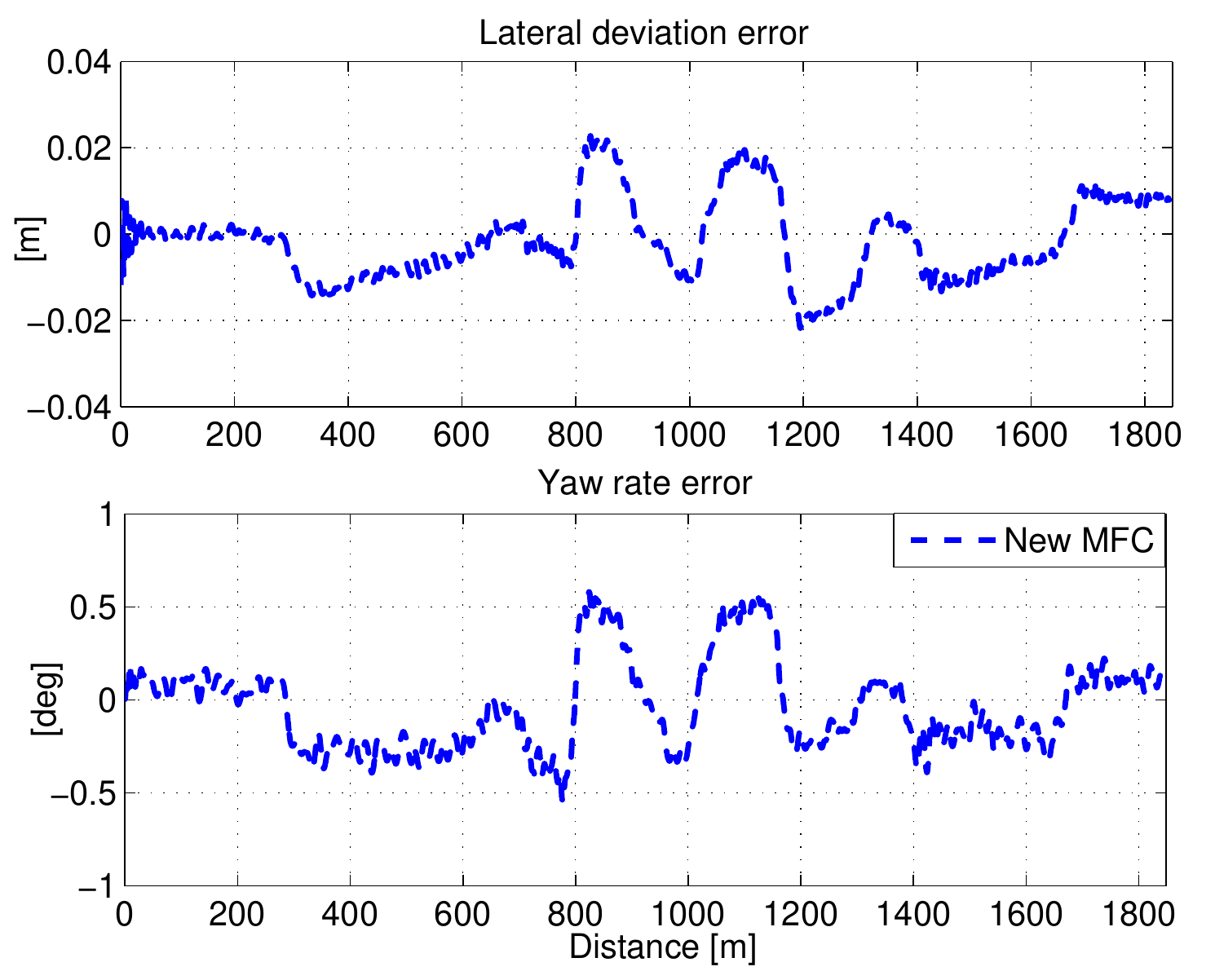}
\vspace{-0.4cm}
\caption{Tracking trajectory errors on lateral deviation and yaw angle}
\label{Ey_Epsi}
\end{figure}

\begin{figure}[!ht]
\centering
\includegraphics[scale=0.56]{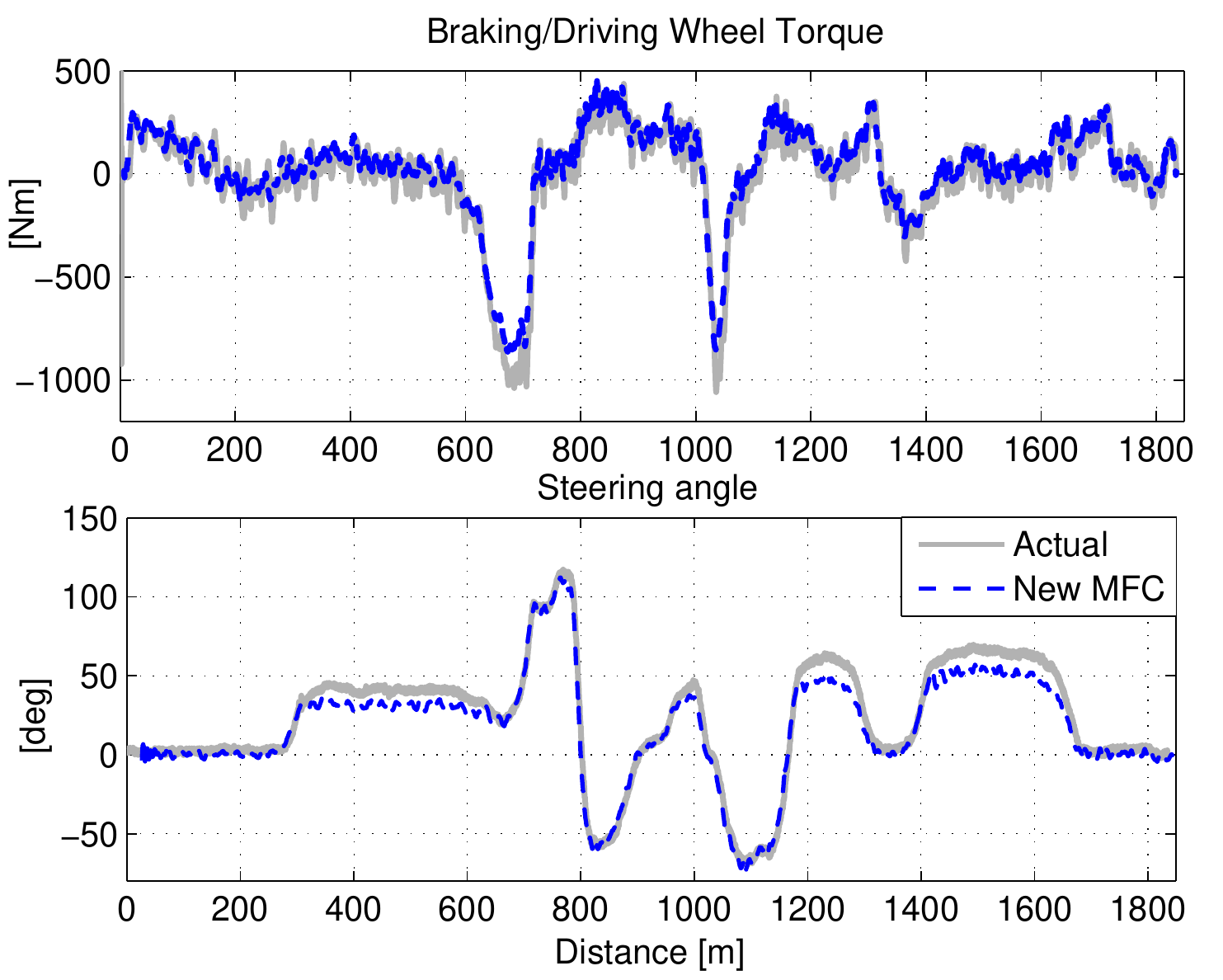}
\vspace{-0.6cm}
\caption{Wheel torques and steering angles control signals: actual and those obtained with our model-free control}
\label{Twheel_Delta}
\end{figure}

Figures \ref{X_Y}, \ref{Ey_Epsi} and \ref{Twheel_Delta} demonstrate that our new model-free setting gives quite satisfying results, even better than the ones obtained in \cite{Menhour13b} with flat outputs. It should be noticed that the test track which has been considered implies strong lateral and longitudinal dynamical requests. This track involves different types of curvatures linked to straight parts, and all these configurations represent a large set of driving situations. Figure \ref{X_Y} shows  that the MFC (model-free control) produces accurate enough behavior for autonomous driving applications. According to the results displayed on Figure \ref{Ey_Epsi}, the lateral error is less than $2 \, cm$. Concerning the yaw angle output, the resulting error is limited to $0.5\, deg$. Finally, Figure \ref{Twheel_Delta} shows that the control signals computed from the MFC strategy are quite closed to the actual ones provided by the driver along the track.

\subsection{Simulation under SiVIC interfaced with RTMaps}

\subsubsection{The SiVIC platform}
Many developments aim to improve the safety of road environments through driving assistance systems. These studies generally take into account an ego vehicle perception and the corresponding vehicle maneuvers (e.g. braking and accelerating). However, in many situations an ego perception is no longer sufficient. Additional information is needed to minimize risk and maximize the security of driving. This additional information requires additional resources which are both time-consuming and expensive. It therefore becomes essential to have a simulation environment allowing to prototype and to evaluate extended, enriched and cooperative driving assistance systems in the early stage of the system design. A virtual simulation platform has to integrate models of road environments, virtual embedded sensors (proprioceptive, exteroceptive), sensors on the infrastructure and communicating devices, according to the laws of physics. In the same way, a physics-based model for vehicle dynamics coupled with actuators (steering wheel angle, torques on each wheel) is necessary. SiVIC meets these requirements and is therefore a very efficient tool to develop and prototype a high level autonomous driving system with cooperative and extended environment perception. 

\begin{figure}[!ht]
\centering
\includegraphics[scale=0.55]{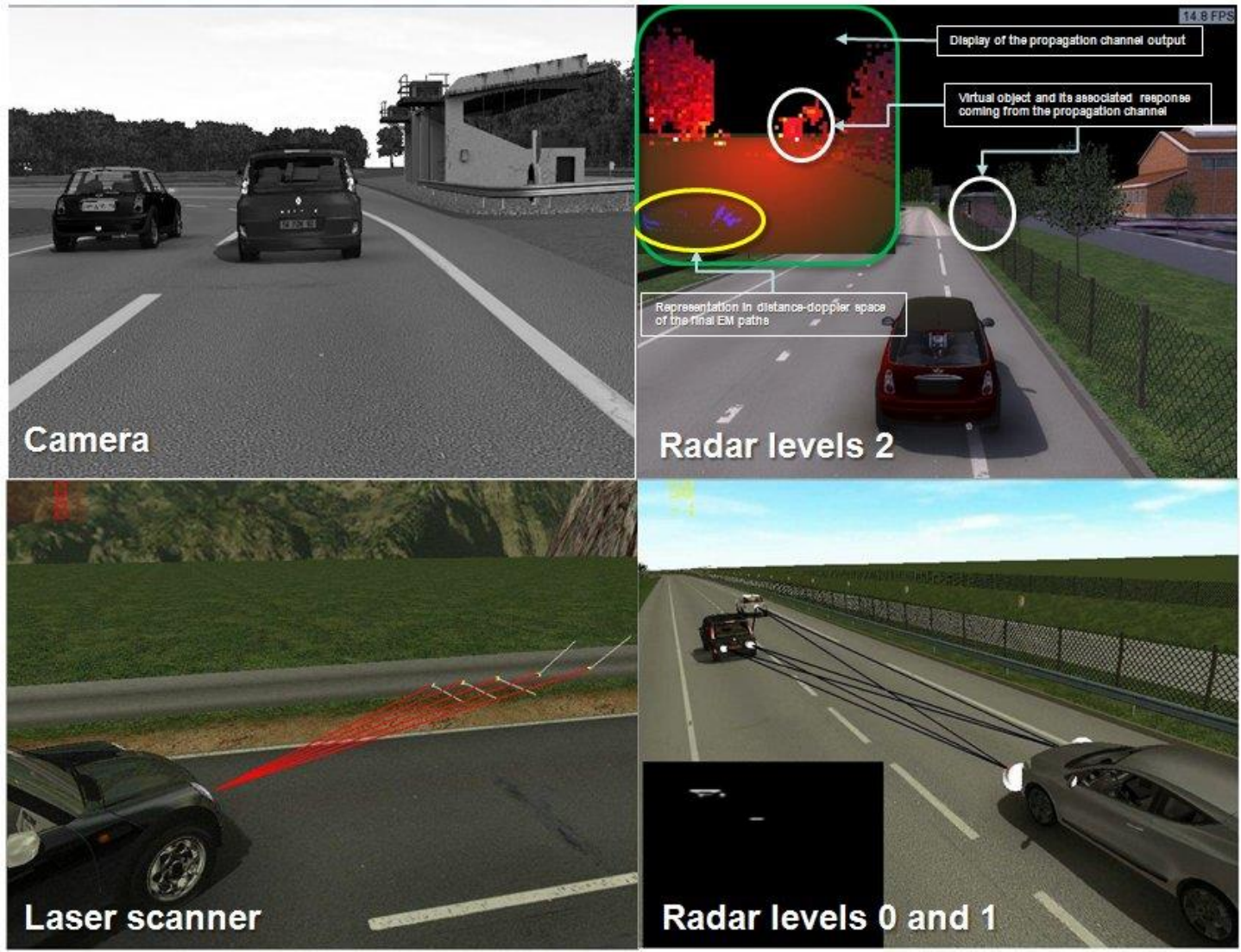}
\vspace{-0.2cm}
\caption{Some exteroceptive types of sensors in SiVIC}
\label{SiVIC_Sensors}
\end{figure}

In its current state, this platform includes several types of exteroceptive and proprioceptive sensors, thus communication media. The exteroceptive sensors are mainly the cameras, the laser scanner, and the RADAR. The proprioceptive sensors involve odometers and Inertial Navigation systems. Then communication sources for cooperative systems include both 802.11p communication media and beacon (transponder). For all this sensors and medias, it is possible in real time and during the simulation stage to tune and to fix the sampling frequency and the intrinsic and extrinsic parameters. Moreover, several modes of operating are available and can be modified during the simulation: ``Off'' and ``On''  to switch on or switch off a sensor. ``Record'' in order to collect data in a file, ``RTMaps'', ``DDS'', and ``Matlab'' to send sensor data in external applications. Some examples of the rendering of these sensors are shown in Fig. \ref{SiVIC_Sensors}.

\subsubsection{RTMaps platform}

RTMaps is a software which has been originally developed at Mines ParisTech by B. Steux \cite{Steux01}.\footnote{It is now edited by the company Intempora.} This platform has been developed for the real-time, multi-sensors, advanced prototyping. Its main goal is to manage and to process a great number of raw data flows like images, laser scanner, GPS, odometric, and INS raw data. The algorithms, which can be applied to the sensor data, are involved in several dedicated image processing and multi-sensors fusion libraries (RTMaps packages). Once these data are recorded and processed, it is very easy to replay them. This type of architecture gives a powerful tool in order to prototype embedded ADAS with either informative outputs or orders to control vehicle dynamics. At each stage, the sensor data and module outputs are time-stamped for an accurate and a reliable time management.

\subsubsection{SiVIC/RTMaps: an interconnected platform for efficient ADAS prototyping}

\begin{figure}[!ht]
\centering
\includegraphics[scale=0.34]{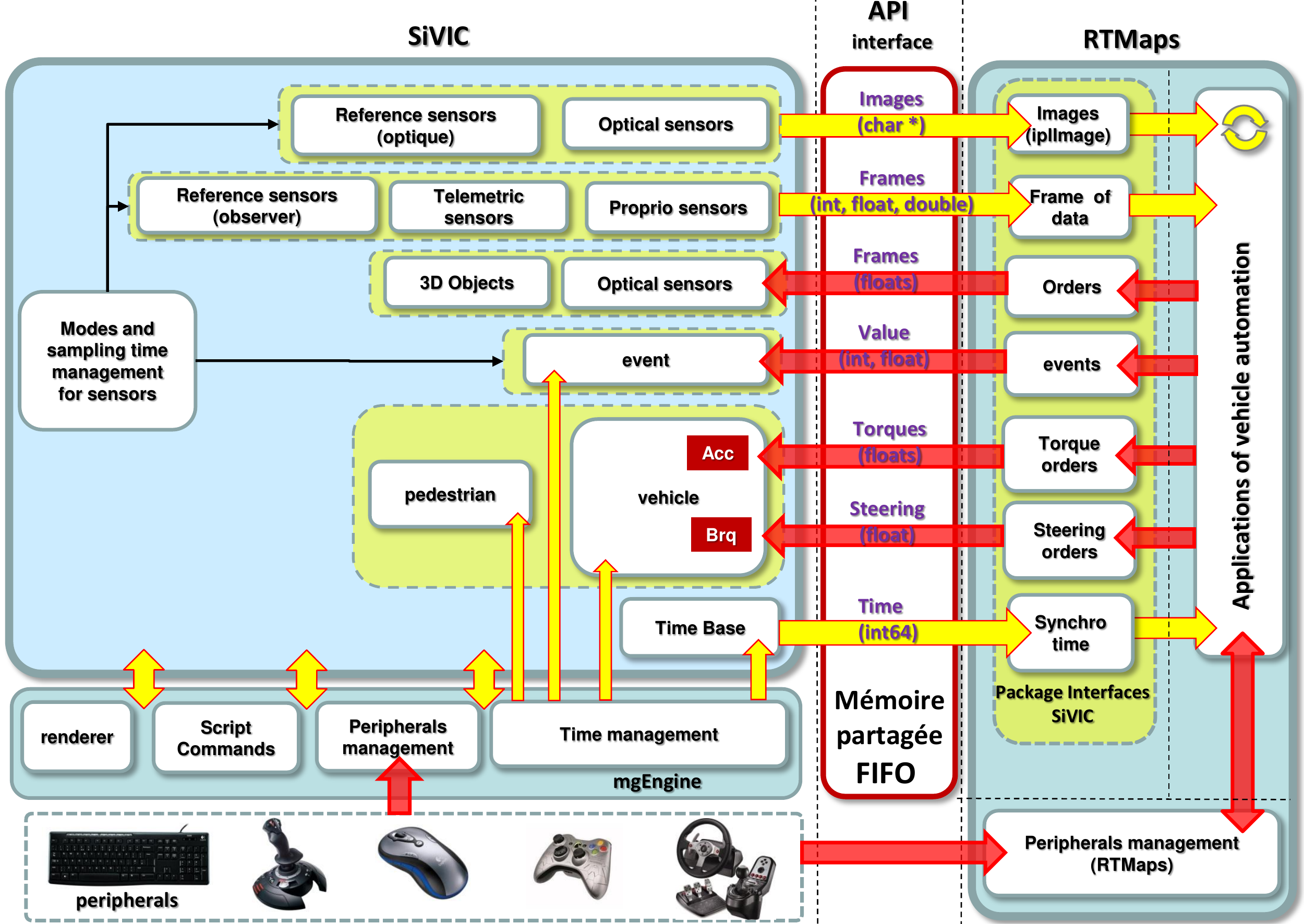}
\vspace{-0.2cm}
\caption{SiVIC/Types of data managed between SiVIC and RTMaps}
\label{SiVIC_Interconnections}
\end{figure}

The interconnection of SiVIC with RTMaps brings RTMaps the ability to replace real-life data by simulated data. Moreover, these interconnected platforms provide a solid framework for advanced prototyping and validation of the control/command and perception algorithms. Indeed, this coupling fully and very effectively allows developing SIL applications (Software In the Loop) including virtual prototypes of vehicles with their proprioceptive and exteroceptive embedded sensors.
The real-time virtual data coming from vehicles and sensors modeled in SiVIC are sent to RTMaps. In RTMaps platform, these data can be used as inputs for perception algorithms and control/command modules. Similarly, orders can be sent from RTMaps to virtual vehicles used in SiVIC in order to control them. This chain of design is very efficient because the algorithms developed in RTMaps can then be directly transferred as micro-software on real hardware devices. Therefore, the simulation model can be considered very close to reality (real vehicles, real sensors). The different types of data handled by this interconnection library are shown in Fig. \ref{SiVIC_Interconnections}.

Several mechanisms have been implemented and tested. The best solution is clearly the optimized FIFO method which allows the transfer of a great number of data in a short time. It is a very critical functionality in order to guarantee a real-time link between SiVIC and the perception/data processing/control algorithms. In order to correctly manage time, a synchronization module is available. This synchronization allows providing a time reference from SiVIC to RTMaps. Then RTMaps is fully synchronized with SiVIC components (vehicle, pedestrian and sensors). The SiVIC/RTMaps simulation platform also enables to build reference scenarios and allows evaluating and testing of control/command and perception algorithms. In fact, the SiVIC/RTMaps platform constitutes a full simulation environment because it provides the same types of interactivity found on actual vehicles: steering wheel angle, acceleration/braking torques, etc.

From the different modules and functionalities available in both SiVIC and RTMaps platform, we have implemented a complete operational architecture in order to test and to evaluate the Model-Free controller with a real time generated reference. In this architecture, shown in Fig. 6 and Fig. 7, SiVIC provides the environment and complex vehicle modeling. Then data coming from a vehicle observer are sent towards RTMaps platform by using a dedicated interface library. In RTMaps, we take into account these car observer data in order to calculate the longitudinal and lateral control inputs. Then, these inputs are sent back to SiVIC's vehicule and more specifically to the virtual actuators. In this architecture, it is also possible to add an event management to apply some constraints in the simulation (obstacle appearance, vehicle parameter modification, ...), and road side beacons to provide speed limit for different areas. The real time software implementation is presented in Fig. 6.

\begin{figure}[!ht]
\centering
\includegraphics[scale=0.35]{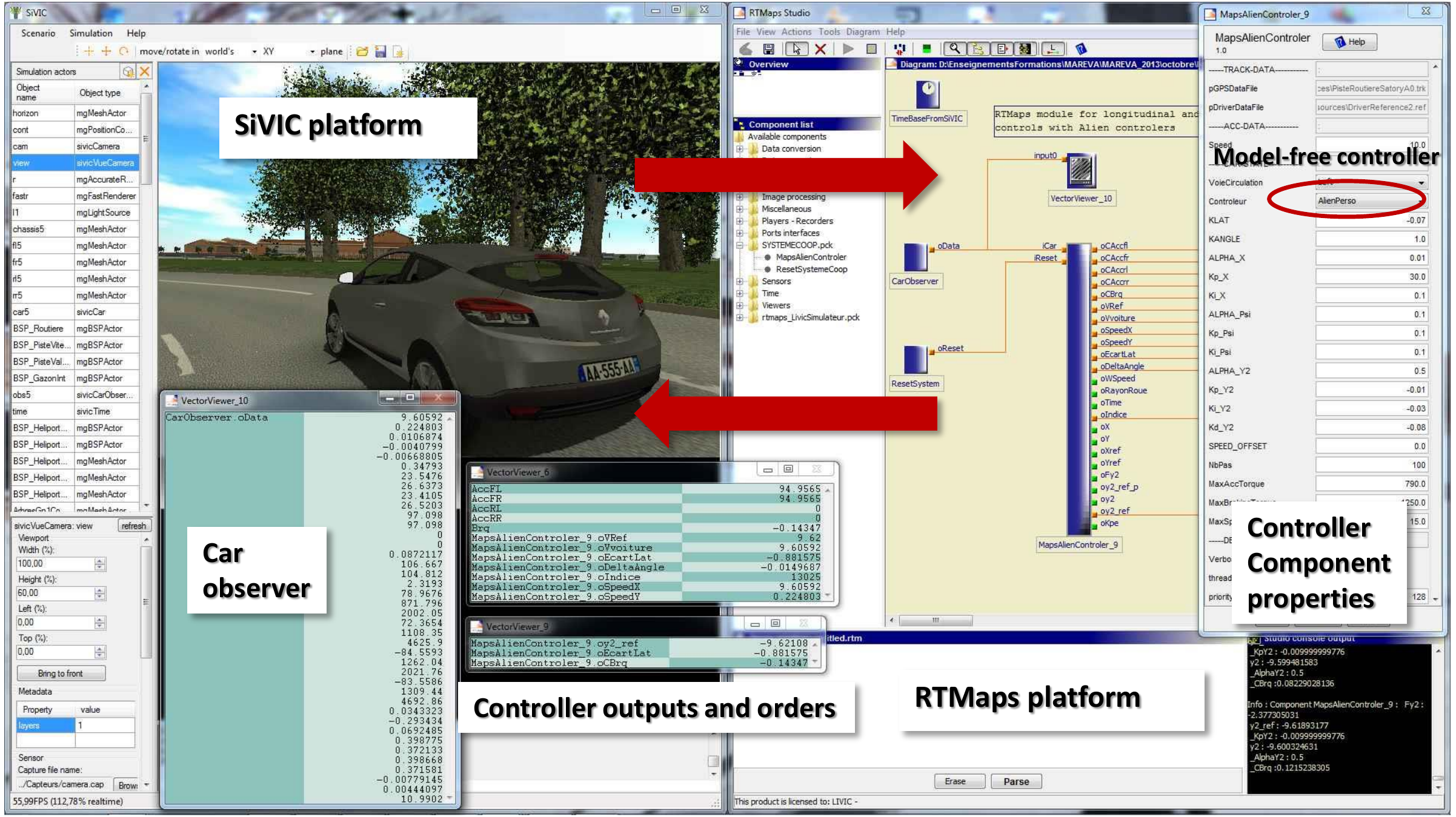}
\vspace{-0.2cm}
\caption{Real time implementation of the model-free control with interconnected platforms RTMAPS/SiVIC}
\label{Inter_Sivic_RTMAPS}
\end{figure}

\begin{figure}[!ht]
\centering
\includegraphics[scale=0.44]{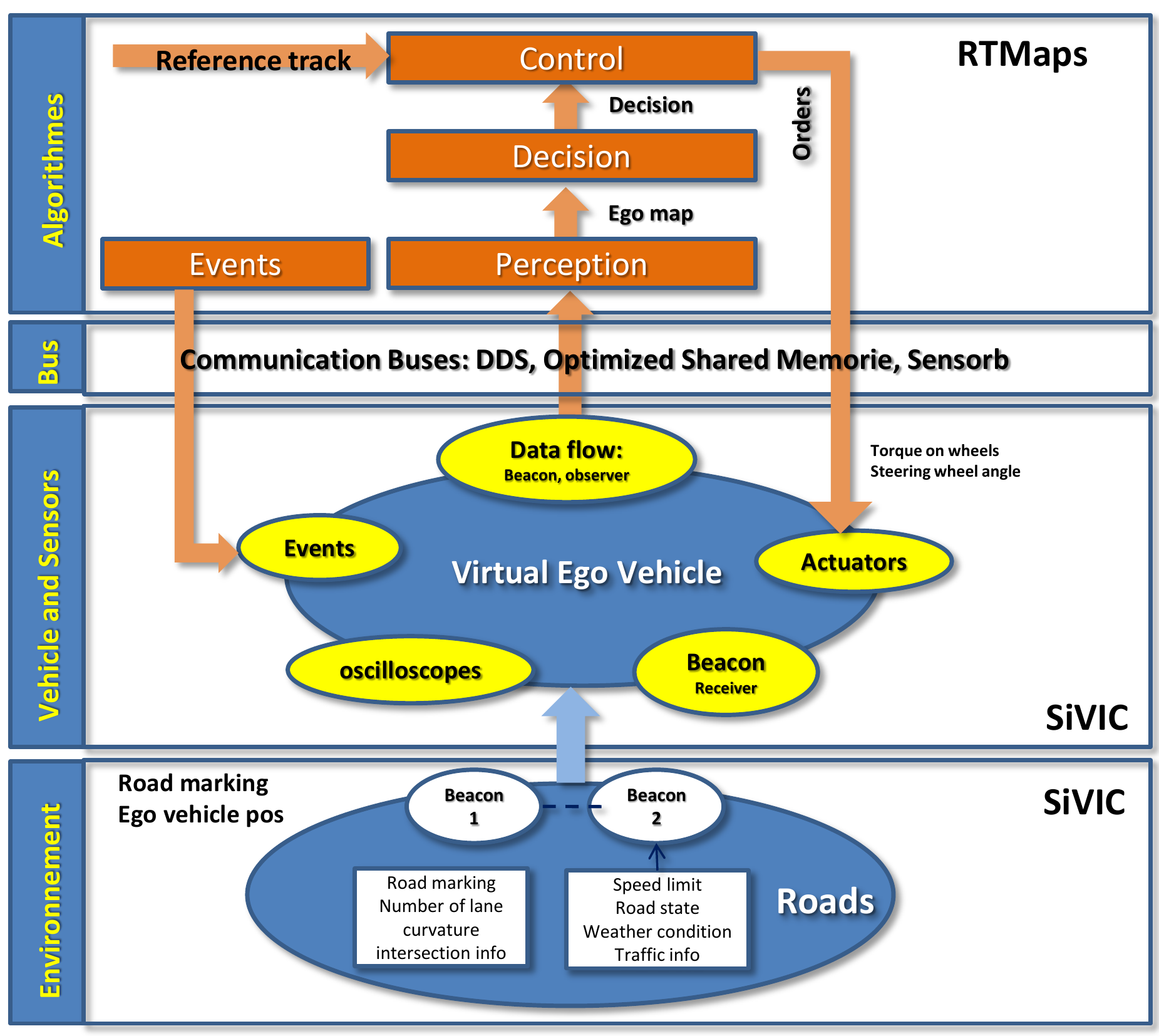}
\vspace{-0.2cm}
\caption{General diagram of software in the loop model-free controller with functional layers}
\label{Inter_Sivic}
\end{figure}

\subsubsection{Simulation results with the interconnected platforms SiVIC and RTMAPS}
\begin{figure}[!ht]
\centering
\includegraphics[scale=0.57]{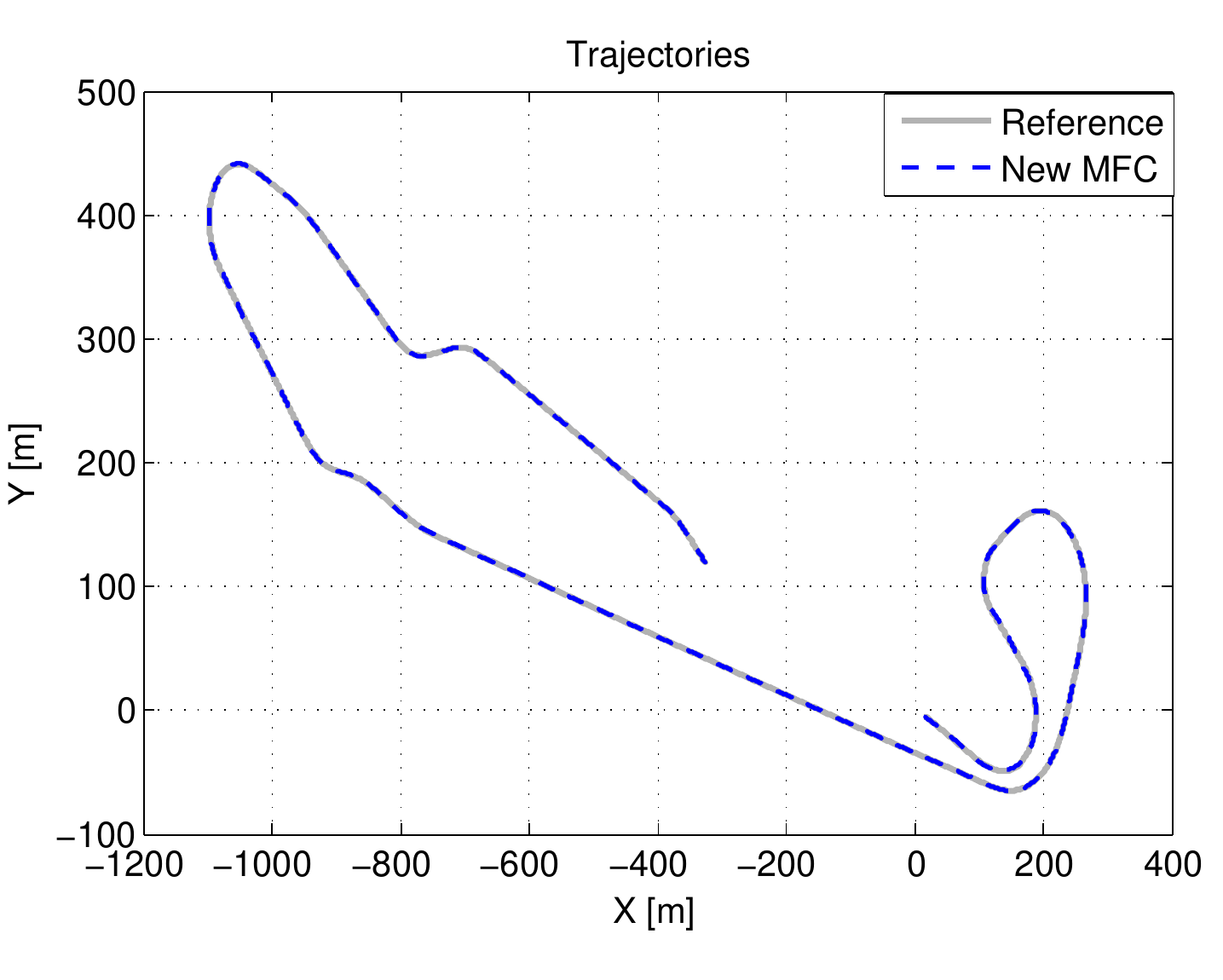}
\caption{Reference trajectory versus the simulated closed-loop trajectory}
\label{X_Y_Sivic}
\end{figure}

\begin{figure}[!ht]
\centering
\includegraphics[scale=0.57]{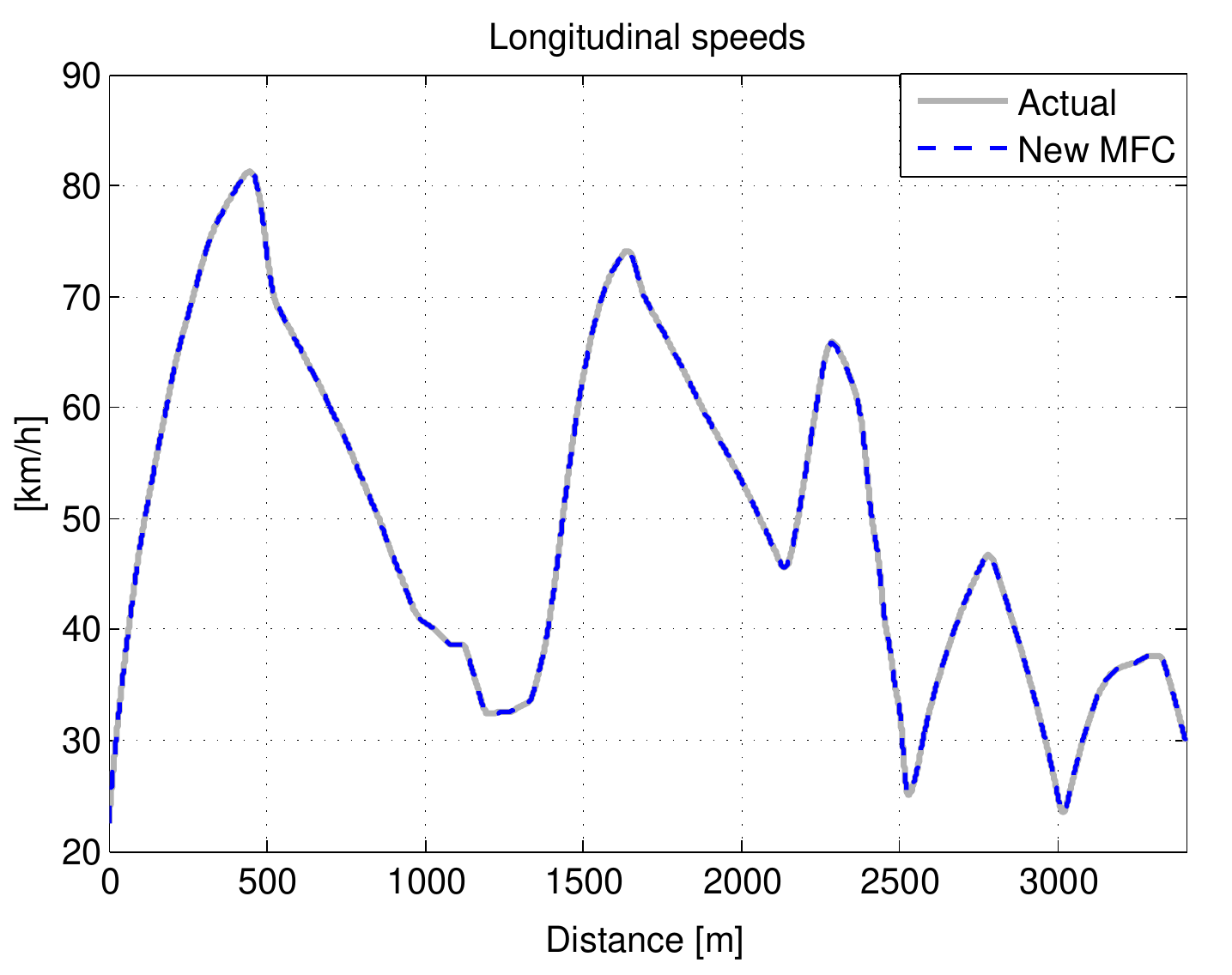}
\caption{Longitudinal speed: actual and simulated}
\label{Vx_SiVIC}
\end{figure}

\begin{figure}[!ht]
\centering
\includegraphics[scale=0.57]{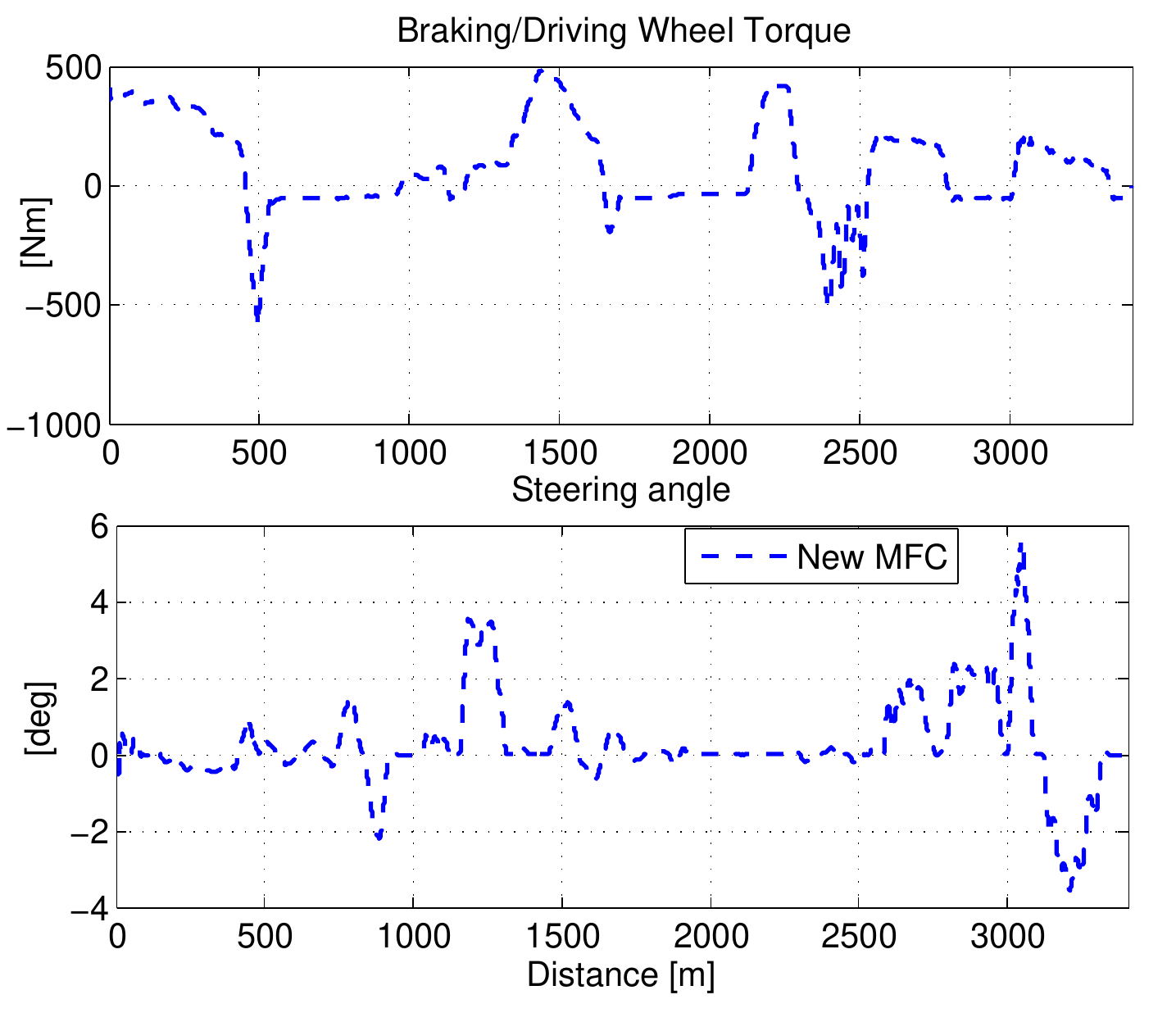}
\caption{Control inputs: Wheel torques and steering angles control signals}
\label{Tw_Delta_SiVIC}
\end{figure}

\begin{figure}[!ht]
\centering
\includegraphics[scale=0.57]{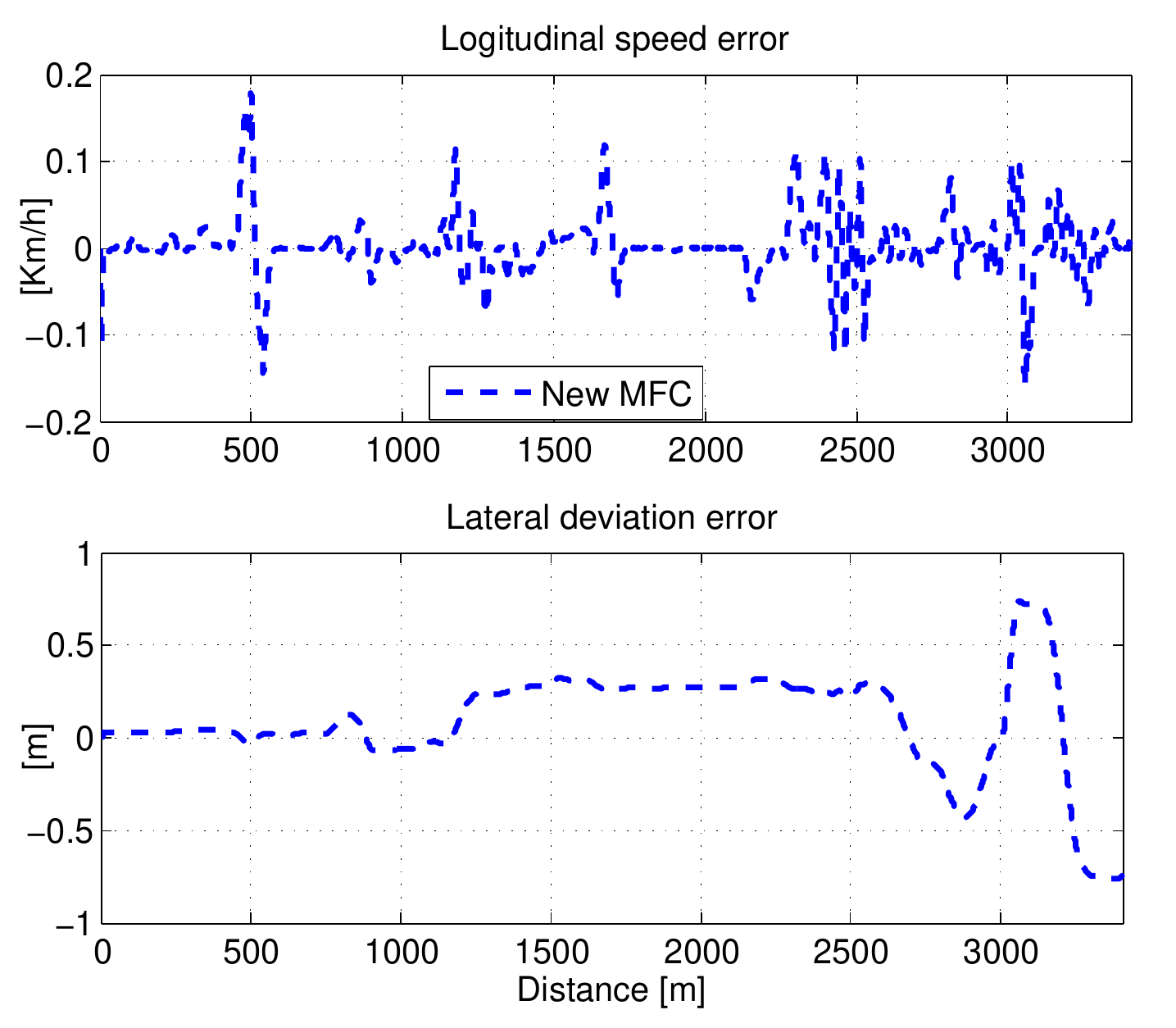}
\caption{Tracking errors on longitudinal speed and lateral deviation}
\label{EVx_Ey_SiVIC}
\end{figure}

Figs. \ref{X_Y_Sivic}, \ref{Vx_SiVIC}, \ref{Tw_Delta_SiVIC} and \ref{EVx_Ey_SiVIC} give the simulation results obtained by implementing the control law under the interconnected SiVIC/RTMaps platform. These results confirm the efficiency of the proposed control law even under a complex and a full simulation environment. The tracking performances, in terms of longitudinal speed and lateral deviation tracking errors, are also depicted in Fig. \ref{EVx_Ey_SiVIC}. With this scenario included very strong curvatures, we can observe a little degradation of the lateral deviation accuracy. Nevertheless the results stay in an acceptable domain allowing to control position of the vehicle in the current lane. About the speed profile, the vehicle follows closely the reference with an absolute error lower than 0.2 km/h. 

\section{Conclusion}\label{Section_6}

For the last past decade, the development of autonomous driving applications is become a key issue for the road transportation domain. In order to solve this problem, a lot of studies and works have been proposed. A great part of these works request an accurate reference and a complex modeling of the system evolution. Unfortunately, the use of a complex model needs the knowledge of both a large quantity of parameters and their operating domain. Often, obtaining these values requests the implementation of expensive experimental means. In this paper, we propose a new model-free approach which takes into account only the outputs signals in order to assess a control law. In application on two different data bases (actual and virtual), this new model-free controller shows that important questions about vehicle control could be efficiently processed with rather elementary and simple tools, which do not necessitate furthermore any complex modeling. In the first data base, actual reference data coming from a Peugeot 406 car has been applied and provide accurate enough results in order to be useful in autonomous driving applications. A second virtual data base has been generated with the SiVIC/RTMaps platform. In the SiVIC platform, an accurate reproduction of the Satory's test tracks is available. Moreover, a complex modeling of the car dynamics is provided. We have interconnected this virtual environment with RTMaps platform in order to obtain a SiL (Software in the Loop) architecture allowing to validate our model-free approach. In this second case, even with very strong curvatures, the results are relevant and accurate.

In future works, we would like to investigate the impact of sensor and observer failures on this method. With this new study, we will confirm that the model-free control might also be robust enough in different types of troubles and noises (see \cite{ijc13,serre}).



\end{document}